\documentclass[11pt]{amsart}
\usepackage{amssymb,amsmath,mathtools}
\usepackage[margin=2.5cm]{geometry}
\usepackage{pdfpages}
\usepackage[colorlinks=true,linkcolor=blue,citecolor=blue]{hyperref}

\title[Calculations for a general maximum principle: notes from 1993]{Calculations for a general maximum principle:\\ 1993 notes on the normalized parabolic $p$-Laplacian}

\author{Nicola Garofalo}
\address{School of Mathematical and Statistical Sciences\\ Arizona State University}
\email{nicola.garofalo@asu.edu}

\date{August 2026}

\begin{document}

\begin{abstract}
This note makes available, in facsimile form, a set of handwritten notes from September 1993 devoted to Bochner-type maximum principle computations for a general class of quasilinear parabolic equations which contains, as a special instance, the normalized parabolic $p$-Laplacian
\[
u_t = \Delta u + (p-2)\,|\nabla u|^{-2}\,\nabla^2u(\nabla u,\nabla u).
\]
The notes were acknowledged in the introduction of the paper of Jin and Silvestre \cite{JS}, and the scheme they propose has recently been implemented in \cite{BGM} to obtain the sharp Li-Yau inequality for this equation. They are reproduced here verbatim, with no alterations, preceded by a brief introduction explaining their content and the subsequent history of the equation.
\end{abstract}

\maketitle

\section{Introduction}

The purpose of this note is to make available a set of handwritten notes which I wrote in September 1993, under the title \emph{Calculations for a general maximum principle}, and which have circulated privately since. The notes are reproduced verbatim, in facsimile, in Section \ref{S:notes} below; nothing in them has been altered, and the notation is that of the time. The present introduction, written in August 2026, has the sole purpose of explaining what the notes contain and why it now seems appropriate to make them publicly available.

\subsection{Content of the notes}
The notes are concerned with the general class of quasilinear parabolic equations
\begin{equation}\label{gen}
a_{ij}(\nabla u)\, \nabla_{ij} u = f(u) + \Phi'(|\nabla u|^2)\, u_t,
\qquad
a_{ij}(\nabla u) = 2\,\Phi''(|\nabla u|^2)\, \nabla_i u\, \nabla_j u + \Phi'(|\nabla u|^2)\,\delta_{ij},
\end{equation}
where $\Phi$ is a given scalar nonlinearity and $f$ a zeroth-order term. The class \eqref{gen} contains the parabolic $p$-Laplacian equation
\begin{equation}\label{me}
|\nabla u|^{p-2}\, u_t = \operatorname{div}\big(|\nabla u|^{p-2}\,\nabla u\big)
\end{equation}
for the choice $\Phi'(s) = s^{\frac{p-2}{2}}$, $f = 0$. Two ideas organize the computations. The first is the normalization of the second-order coefficients by the structure function
\[
\Lambda(s) := 2s\,\Phi''(s) + \Phi'(s),
\]
which is the largest eigenvalue of the coefficient matrix $a_{ij}$: dividing \eqref{gen} through by $\Phi'(|\nabla u|^2)$ produces, in the case of \eqref{me}, precisely the \emph{normalized} equation
\begin{equation}\label{norm}
u_t = \Delta u + (p-2)\,|\nabla u|^{-2}\,\nabla^2u(\nabla u,\nabla u),
\end{equation}
whose coefficient matrix $I + (p-2)\,\xi\otimes\xi$, $\xi = \nabla u/|\nabla u|$, has eigenvalues $1$ and $p-1$: unlike \eqref{me}, the equation \eqref{norm} is uniformly elliptic in its dependence on the Hessian, and it was observed in the notes that, for this reason, it should enjoy good regularization properties despite the degeneracy of the divergence form \eqref{me}. The second idea is a maximum principle scheme run on functionals of the form
\[
P = \psi(|\nabla u|^2) - 2F(u),
\]
in which the scalar nonlinearity $\psi$ is \emph{adapted to the operator} (in the notes, $\psi' = \Lambda$): the choice is dictated by the requirement that, after the first-order conditions $\nabla P = 0$, $P_t = 0$ are enforced at an interior maximum point, the dangerous quadratic terms generated by the Bochner-type differentiation cancel exactly. The second-order coefficients and the scalar structure functions are treated throughout as decoupled objects, to be matched at the end of the computation.

\subsection{Subsequent history}
The normalized equation \eqref{norm} has since acquired an independent life. It arises as the evolution governing the value functions of the tug-of-war games with noise of Peres and Sheffield \cite{PS}, and it has been studied by several authors: we refer to \cite{Do, BG, JS} for well-posedness and regularity in the viscosity framework, and to \cite{MPR} for the parabolic game-theoretic interpretation. The notes reproduced below are acknowledged in the introduction of the paper of Jin and Silvestre \cite{JS}, where it is noted that they contain a computation leading to Lemma 3.1 of that paper, and that they constitute, to the authors' knowledge, the first appearance of the equation \eqref{norm} together with the recognition of its regularization properties.

\subsection{Why now}
In the recent work \cite{BGM} with A.~Banerjee and H.~Mahmoudian we prove the sharp Li-Yau inequality
\[
(p-1)\frac{|\nabla u|^2}{u^2} - \frac{u_t}{u} \;\le\; \frac{n+p-2}{2(p-1)\,t}
\]
for positive viscosity solutions of \eqref{norm} on closed Riemannian manifolds with nonnegative Ricci curvature, with no assumption on the critical set of the solution. The proof rests on a family of uniformly parabolic approximating flows in which the second-order regularization and the first-order term are decoupled, the latter being chosen so that the dangerous terms in the Bochner identity cancel exactly and the sharp constant survives the approximation. This is, in the Li-Yau setting, an implementation of the scheme proposed in the 1993 notes, and \cite{BGM} refers to them at several points. Since the notes have never been publicly available, it seems appropriate -- for completeness of the record, and because the general class \eqref{gen} may still hold interest beyond the special case \eqref{norm} -- to make them accessible in their original form. This is done in the next section.

\newpage

\section{The notes}\label{S:notes}

The following ten pages reproduce, in facsimile, the handwritten notes \emph{Calculations for a general maximum principle}, September 1993. A word on the facsimile is in order. The pages are presented here in their original reading order, which differs from the order in which they were scanned. They consist of three parts: a preliminary check in one space dimension (first page); the main computation, opening with the words \emph{``We start with the normalized parabolic equation''} and comprising equations (1)--(16) together with its unnumbered conclusion (six pages); and a second, refined computation, dated September 6, 1993, comprising equations (31)--(41) and its conclusion (four pages). The jump in the numbering, from (16) to (31), indicates that some intermediate pages have not survived the intervening years; the two computations reproduced here are, however, each self-contained. Nothing has been altered.

\includepdf[pages={5-10,1-4},scale=0.95,pagecommand={}]{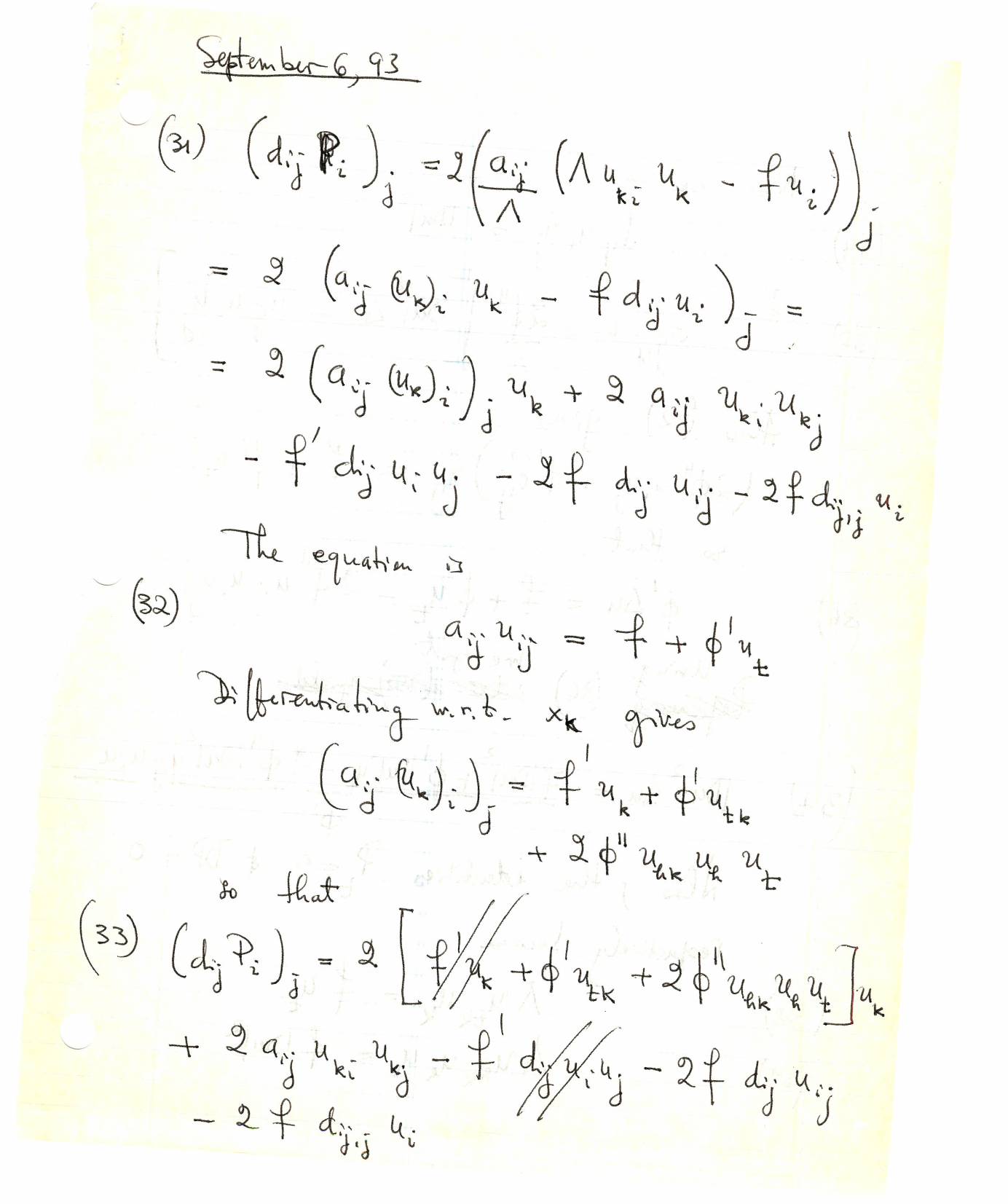}

\end{document}